\documentclass{amsart}
\usepackage{amsmath,amssymb,amsthm}
\usepackage{bm}
\usepackage{enumerate}
\usepackage{extarrows}
\usepackage{mathrsfs}
\usepackage{stmaryrd}
\usepackage[colorlinks=true]{hyperref}
\usepackage[all]{xy}

\newtheorem{theorem}{Theorem}[section]
\newtheorem{corollary}[theorem]{Corollary}
\newtheorem{lemma}[theorem]{Lemma}
\newtheorem{proposition}[theorem]{Proposition}

\theoremstyle{definition}

\newtheorem{example}[theorem]{Example}

\theoremstyle{remark}

\numberwithin{equation}{section}

\newcommand\ra{\rightarrow}
\newcommand\lra{\longrightarrow}

\newcommand\simto{\sto{\sim}}
\newcommand\lto{\longmapsto}

\font\cyr=wncyr10
\newcommand\Sha{\hbox{\cyr X}}

\newcommand\sto[1]{\stackrel{#1}{\longrightarrow}}

\newcommand\fct[4]{\begin{split}#1 &\lra #2 \\ #3 &\lto #4\end{split}}
\newcommand\set[1]{{\left\{#1\right\}}}
\newcommand\pair[1]{\langle{#1}\rangle}

\newcommand\leg[2]{\left(\frac{{#1}}{#2}\right)}
\newcommand\aleg[2]{\left[\frac{{#1}}{#2}\right]}

\newcommand\rmT{{\mathrm{T}}}

\newcommand\bfb{{\mathbf{b}}}

\newcommand\bfd{{\mathbf{d}}}

\newcommand\bfx{{\mathbf{x}}}
\newcommand\bfy{{\mathbf{y}}}

\newcommand\bfA{{\mathbf{A}}}
\newcommand\bfB{{\mathbf{B}}}

\newcommand\bfD{{\mathbf{D}}}

\newcommand\bfM{{\mathbf{M}}}
\newcommand\bfN{{\mathbf{N}}}
\newcommand\bfO{{\mathbf{O}}}

\newcommand\bfR{{\mathbf{R}}}

\newcommand\BA{{\mathbb{A}}}

\newcommand\BF{{\mathbb{F}}}

\newcommand\BN{{\mathbb{N}}}

\newcommand\BQ{{\mathbb{Q}}}

\newcommand\BZ{{\mathbb{Z}}}

\newcommand\CO{{\mathcal{O}}}

\newcommand\fd{{\mathfrak{d}}}

\newcommand\cl{{\mathrm{cl}}}

\newcommand\diag{{\mathrm{diag}}}

\renewcommand\Im{{\mathrm{Im}\,}}

\newcommand\Ker{{\mathrm{Ker}\,}}

\renewcommand\mod{\, \mathrm{mod}\, }

\newcommand\rank{{\mathrm{rank}}}

\newcommand\Sel{{\mathrm{Sel}}}

\title[On non-congruent numbers and tame kernels]{On non-congruent numbers with $8a\pm1$ \\ type odd prime factors and tame kernels}
\author{Shenxing Zhang}
\date{\today}
\address{School of Mathematics, Hefei University of Technology, Hefei, Anhui 230000, China}
\email{zhangshenxing@hfut.edu.cn}

\begin{document}

\begin{abstract}
Let $n$ be a positive square-free integer, where every odd prime factor of $n$ has form $8a\pm 1$.
We determine when $n$ is non-congruent with second minimal $2$-primary Shafarevich-Tate group, in terms of the $4$-ranks of class groups and a Jacobi symbol.
In particular, when every odd prime factor of $n$ has form $8a+1$, this condition is equivalent to the vanishing of the $4$-rank of the tame kernel of $\BQ(\sqrt{n})$ for odd $n$, or $\BQ(\sqrt{-n})$ for even $n$. This generalizes previous results.
\end{abstract}

\maketitle
\tableofcontents

\section{Introduction}

\subsection{Background}
A square-free positive integer $n$ is called \emph{congruent} if it is the area of a right triangle with rational lengths.
This is equivalent to say, the Mordell-Weil rank of $E$ over $\BQ$ is positive, where
\begin{equation}\label{eq:congruent_ec}
	E=E_n: y^2=x^3-n^2x
\end{equation}
is the associated congruent elliptic curve.
Denote by $\Sel_2(E)$ the \emph{$2$-Selmer group} of $E$ over $\BQ$ and 
\begin{equation}\label{eq:pure2selmer_rank}
	s_2(n):=\dim_{\BF_2}\left(\frac{\Sel_2(E)}{E(\BQ)[2]}\right)
		=\dim_{\BF_2}\Sel_2(E)-2
\end{equation}
the \emph{pure $2$-Selmer rank}.
Then
\[s_2(n)=\rank_\BZ E(\BQ)+\dim_{\BF_2}\Sha(E/\BQ)[2]\]
by the exact sequence
\[0\ra E(\BQ)/2E(\BQ)\ra \Sel_2(E)\ra \Sha(E/\BQ)[2]\ra0.\]

Certainly, $s_2(n)=0$ implies that $n$ is non-congruent.
The examples of $s_2(n)=0$ can be found in \cite{Feng1997}, \cite{Iskra1996} and \cite{OuyangZhang2015}, which are corollaries of Monsky's formula for $s_2(n)$. This case is fully characterized in terms of class groups and the full Birch-Swinnerton-Dyer (BSD) conjecture holds, see \cite[Theorem~1.1, Corollary~1.3]{TianYuanZhang2017} and \cite[Theorem~1.2]{Smith2016}.
The examples of non-congruent $n$ with $\Sha(E/\BQ)[2^\infty]\cong(\BZ/2\BZ)^2$ can be found in \cite{LiTian2000}, \cite{OuyangZhang2014} and \cite{OuyangZhang2015}.
When $n$ is odd with prime factors $1\bmod 4$, it can be characterized as follows.

Denote by $(a,b)_v$ the Hilbert symbol.
Denote by 
\begin{equation}\label{eq:2rank}
	r_{2^a}(A)=\dim_{\BF_2}\left(\frac{2^{a-1}A}{2^aA}\right)
\end{equation}
the $2^a$-rank of a finite abelian group $A$.
Denote by $h_{2^a}(m)$ the $2^a$-rank of the narrow class group of $\BQ(\sqrt{m})$.

\begin{theorem}[{\cite[Theorem~1.1]{Wang2016}}]
Let $n=p_1\cdots p_k\equiv 1\bmod 8$ be a square-free positive integer with $p_i\equiv 1\bmod 4$.
The following are equivalent:
\begin{enumerate}[(i)]
\item $h_4(-n)=1$ and $h_8(-n)\equiv \frac{d-1}{4}\pmod 2$;
\item $\rank_\BZ E_n(\BQ)=0$ and $\Sha(E_n/\BQ)[2^\infty]\cong(\BZ/2\BZ)^2.$
\end{enumerate}
Here, either $d\neq 1,n$ is a positive divisor of $n$ such that $(d,-n)_v=1,\forall v$, or $d$ is a positive divisor of $n$ such that $(2d,-n)_v=1,\forall v$.
\end{theorem}

Wang also gave a sufficient condition on $\Sha(E/\BQ)\cong(\BZ/2\BZ)^{s_2(n)}$ for $s_2(n)\ge 4$.
Recently, Qin in \cite[Theorem~1.5]{Qin2021} proved that if $p\equiv 1\bmod 8$ is a prime with trivial $8$-rank of the tame kernel $K_2\CO_{\BQ(\sqrt{p})}$, then $p$ is non-congruent.
Moreover, if the $4$-rank of $K_2\CO_{\BQ(\sqrt{p})}$ is $1$, then $\Sha(E_p/\BQ)[2^\infty]\cong(\BZ/4\BZ)^2$.

\subsection{Main results}
In this paper, we will give a criterion of non-congruent $n$ with $\Sha(E_n/\BQ)[2^\infty]\cong(\BZ/2\BZ)^2$, where the odd prime factors of $n$ are congruent to $\pm 1$ modulo $8$.

\begin{theorem}[=Theorem~\ref{thm:main1}]
Let $n=p_1\cdots p_k\equiv 1\bmod 8$ be a square-free positive integer with $p_i\equiv \pm1\bmod8$.
The following are equivalent:
\begin{enumerate}[(i)]
\item $h_4(-n)=1, h_4(n)=0$ and $\leg{-\mu}{d}=-1$;
\item $\rank_\BZ E_n(\BQ)=0$ and $\Sha(E_n/\BQ)[2^\infty]\cong(\BZ/2\BZ)^2$.
\end{enumerate}
Here, $d\neq 1$ is the unique positive divisor of $n$ such that $(d,n)_v=1,\forall v$ and $n=2\mu^2-\tau^2$ where $\mu\equiv d\bmod4$.
\end{theorem}

\begin{corollary}[= Corollary~\ref{cor:main1}]
Let $n=p_1\cdots p_k\equiv 1\bmod 8$ be a square-free positive integer with $p_i\equiv 1\bmod8$.
The following are equivalent:
\begin{enumerate}[(i)]
\item $r_4(K_2\CO_{\BQ(\sqrt{n})})=0$;
\item $\rank_\BZ E_n(\BQ)=0$ and $\Sha(E_n/\BQ)[2^\infty]\cong(\BZ/2\BZ)^2$.
\end{enumerate}
\end{corollary}
This generalizes \cite[Lemma~4.2]{Qin2021}.

\begin{theorem}[= Theorem~\ref{thm:main2}]
Let $n=2p_1\cdots p_k\equiv 2\bmod 8$ be a square-free positive integer with $p_i\equiv \pm1\bmod8$.
The following are equivalent:
\begin{enumerate}[(i)]
\item $h_4(-n/2)=1$ and $\leg{2-\sqrt{2}}{|d|}=-1$;
\item $\rank_\BZ E_n(\BQ)=0$ and $\Sha(E_n/\BQ)[2^\infty]\cong(\BZ/2\BZ)^2$.
\end{enumerate}
Here, $d\neq 1$ is the unique divisor of $n$ such that $(d,n)_v=1,\forall v$ and $d\equiv1\bmod8$.
\end{theorem}

\begin{corollary}[= Corollary~\ref{cor:main2}]
Let $n=2p_1\cdots p_k\equiv 2\bmod 8$ be a square-free positive integer with $p_i\equiv 1\bmod8$.
The following are equivalent:
\begin{enumerate}[(i)]
\item $r_4(K_2\CO_{\BQ(\sqrt{-n})})=0$;
\item $\rank_\BZ E_n(\BQ)=0$ and $\Sha(E_n/\BQ)[2^\infty]\cong(\BZ/2\BZ)^2$.
\end{enumerate}
\end{corollary}

\subsection{Notations}

\begin{itemize}
\item $E=E_n$ the congruent elliptic curve associated a square-free positive integer $n$, see \eqref{eq:congruent_ec}.
\item $s_2(n)$ the pure $2$-Selmer rank of $E_n$, see \eqref{eq:pure2selmer_rank}.
\item $(m,n)$ the greatest common divisor of integers $(m,n)$, where $m\neq0$ or $n\neq 0$.
\item $(a,b)_v$ the Hilbert symbol.
\item $[a,b]_v$ the additive Hilbert symbol, i.e., the image of $(a,b)_v$ under the isomorphism $\set{\pm1}\simto\BF_2$.
\item $\leg{a}{b}=\prod_{p\mid b}(a,b)_p$ the Jacobi symbol, where $(a,b)=1$ and $b>0$.
\item $\aleg{a}{b}$ the additive Jacobi symbol, i.e., the image of $\leg{a}{b}$ under the isomorphism $\set{\pm1}\simto\BF_2$.
\item $r_{2^a}$ the $2^a$-rank of a finite abelian group, see \eqref{eq:2rank}.
\item $h_{2^a}(m)$ the $2^a$-rank of the narrow class group of $\BQ(\sqrt{m})$.
\item $K_2\CO_F$ the tame kernel of a number field $F$, see \cite[Theorem~III.6.5]{Weibel2013}.
\item $m'=|m|/(2^{|m|},m)$ the positive odd part of an integer $m$.
\item $\bfA=\bfA_{n'}$ a matrix associated to $n'$, see \eqref{eq:monsky}.
\item $\bfD_\varepsilon$ a matrix associated to $n'$ and $\varepsilon$, see \eqref{eq:D_epsilon}.
\item  ${\bf0}=(0,\dots,0)^\rmT$, ${\bf1}=(1,\dots,1)^\rmT$ and $\bfb_\varepsilon=\bfD_\varepsilon{\bf1}$.
\item $\bfM_n$ the Monsky matrix of $E_n$, see \eqref{eq:monsky_matrix1} and \eqref{eq:monsky_matrix1}.
\item $\bfR_m$ the R\'edei matrix of $\BQ(\sqrt{m})$, see \eqref{eq:Redei_matrix}.
\item $V,V_1,V_2$ sets associated to $\BQ(\sqrt{m})$, see Theorems~\ref{thm:K_2_odd} and \ref{thm:K_2_even}.
\item $v_p$ the normalized valuation on $\BQ_p$.
\item For a vector $\bfd=(\delta_1,\dots,\delta_k)^\rmT\in\BF_2^k$, denote by $\psi(\bfd)=\prod_{j=1}^k p_j^{\delta_j}$. Then $\psi^{-1}(|d|)=\bigl(v_{p_1}(d),\dots,v_{p_k}(d)\bigr)^\rmT$ for $d\mid p_1\dots p_k$.
\end{itemize}

\section{Preliminaries}

\subsection{Monsky matrix}
\label{ssec:monsky}
Monsky in \cite[Appendix]{HeathBrown1994} represented the pure $2$-Selmer group as the kernel of a matrix over $\BF_2$.
Let's recall it roughly.
One can identify $\Sel_2(E_n)$ with
\[\set{\Lambda=(d_1,d_2,d_3)\in(\BQ^\times/\BQ^{\times2})^3:
D_\Lambda(\BA_\BQ)\neq \emptyset,d_1d_2d_3\equiv 1\bmod\BQ^{\times 2}},\]
where $D_\Lambda$ is a genus one curve defined by
\begin{equation}
	\begin{cases}
		H_1:& -nt^2+d_2u_2^2-d_3u_3^2=0,\\
		H_2:& -nt^2+d_3u_3^2-d_1u_1^2=0,\\
		H_3:& 2nt^2+d_1u_1^2-d_2u_2^2=0.
	\end{cases}
\end{equation}
Under this identification, $O,(n,0),(-n,0),(0,0)$ and non-torsion  $(x,y)\in E_n(\BQ)$ correspond to $(1,1,1),(2,2n,n),(-2n,2,-n),(-n,n,-1)$ and $(x-n,x+n,x)$ respectively.

Let $n'=p_1\cdots p_k$ be an ordered prime decomposition of $n'=n/(2,n)$.
Denote by
\begin{equation}\label{eq:monsky}
	\bfA=\bfA_{n'}=(a_{ij})_{k\times k}
	\quad\text{where}\quad
	a_{ij}=[p_j,-n']_{p_i}=
	\begin{cases}
		\aleg{p_j}{p_i},&i\neq j;\\
		\aleg{n'/p_i}{p_i},&i=j,
	\end{cases}
\end{equation}
and
\begin{equation}\label{eq:D_epsilon}
\bfD_\varepsilon=\diag\left\{\aleg{\varepsilon}{p_1},\dots,\aleg{\varepsilon}{p_k}\right\}.
\end{equation}
Then $\bfA{\bf1}={\bf0}$ and $\rank(\bfA)\le k-1$.

When $n$ is odd, each element in $\Sel_2(E_n)/E_n(\BQ)[2]$ can be presented as $(d_1,d_2,d_3)$, where both of $d_1,d_2$ are positive divisor of $n$.
The system $D_\Lambda$ is locally solvable everywhere if and only if certain conditions on Hilbert symbols hold.
Then we can express the pure $2$-Selmer group as the kernel of Monsky matrix
\begin{equation}\label{eq:monsky_matrix1}
	\bfM_n=\begin{pmatrix}
		\bfA+\bfD_2&\bfD_2\\
		\bfD_2&\bfA+\bfD_{-2}
	\end{pmatrix}
\end{equation}
via the isomorphism
\begin{equation}\label{eq:monsky_isom1}
\fct{\Sel_2(E_n)/E_n(\BQ)[2]}{\Ker\bfM_n}{(d_1,d_2,d_3)}{\left(\begin{smallmatrix}
\psi^{-1}(d_2)\\ \psi^{-1}(d_1)
\end{smallmatrix}\right).}
\end{equation}
% A+A^T=D_{-1}+\bfb_{-1}\bfb_{-1}^\rmT

When $n$ is even, each element in $\Sel_2(E_n)/E_n(\BQ)[2]$ can be presented as $(d_1,d_2,d_3)$, where both of $d_2,d_3$ are divisor of $n'$ and $d_2>0, d_3\equiv1\bmod4$.
Then we can express the pure $2$-Selmer group as the kernel of Monsky matrix
\begin{equation}\label{eq:monsky_matrix2}
	\bfM_n=\begin{pmatrix}
		\bfA^\rmT+\bfD_{2}&\bfD_{-1}\\
		\bfD_2&\bfA+\bfD_2
	\end{pmatrix}
\end{equation}
via the isomorphism
\begin{equation}\label{eq:monsky_isom2}
\fct{\Sel_2(E_n)/E_n(\BQ)[2]}{\Ker\bfM_n}{(d_1,d_2,d_3)}{\left(\begin{smallmatrix}
\psi^{-1}(|d_3|)\\ \psi^{-1}(d_2)
\end{smallmatrix}\right).}
\end{equation}

In both cases, we have
\begin{equation}\label{eq:monsky_rank}
s_2(n)=2k-\rank(\bfM_n).
\end{equation}

\subsection{Cassels pairing}
Cassels in \cite{Cassels1998} defined a skew-symmetric bilinear pairing $\pair{-,-}$ on the $\BF_2$-vector space $\Sel_2(E_n)/E_n(\BQ)[2]$.
For any $\Lambda\in\Sel_2(E_n)$, choose any $P=(P_v)\in D_\Lambda(\BA_\BQ)$.
Note that $H_i$ is locally solvable everywhere, hence it is solvable over $\BQ$ by Hasse-Minkowski principal.
Choose $Q_i\in H_i(\BQ)$. Let $L_i$ be a linear form in three variables such that $L_i=0$ defines the tangent plane of $H_i$ at $Q_i$.
Then for any $\Lambda'=(d'_1,d'_2,d'_3)\in\Sel_2(E_n)$, define
\[\pair{\Lambda,\Lambda'}=\prod_v\pair{\Lambda,\Lambda'}_v
\quad\text{where}\quad \pair{\Lambda,\Lambda'}_v=\prod_{i=1}^3 \bigl(L_i(P_v),d_i'\bigr)_v.\]
This pairing is independent of the choice of $P$ and $Q_i$ and is trivial on $E_n(\BQ)[2]$.

\begin{lemma}[{\cite[Lemma~7.2]{Cassels1998}}]\label{lem:cassels}
The local Cassels pairing $\pair{-,-}_p=+1$ if 
\begin{itemize}
\item $p\nmid 2\infty$,
\item the coefficients of $H_i$ and $L_i$ are all integral at $p$ for $i=1,2,3$, and
\item modulo $D_\Lambda$ and $L_i$ by $p$, they define a curve of genus $1$ over $\BF_p$ together with tangents to it.
\end{itemize}
\end{lemma}

\subsection{R\'edei matrix}
Let $m\neq 0,1$ be a square-free integer.
Denote by $F=\BQ(\sqrt{m})$ and $\bfN F=\bfN_{F/\BQ}(F^\times)$.
Denote by $C(F)$ the narrow class group of $F$.
Let $D=p_1^*\cdots p_t^*$ be the decomposition of the discriminant of $F$ into a product of \emph{prime discriminants}
\[p^*=(-1)^{\frac{p-1}2}p,\qquad 2^*=-4,8,-8.\]
By Gauss genus theory, we have $h_2(m)=t-1$.	
To calculate $h_4(m)$, we need the R\'edei matrix, which is defined as
\begin{equation}\label{eq:Redei_matrix}
	\bfR_m=([p_j,m]_{p_i})_{t\times t}.
\end{equation}
Let $V$ be the set of all square-free positive integers $d\mid D$.
Then the following are equivalent:
\begin{itemize}
\item $d\in V\cap \bfN F$;
\item $X^2-mY^2=dZ^2$ is has nonzero solutions over $\BZ$;
\item Hilbert symbols $(d,m)_p=1, \forall p\mid D$;
\item $\bfR_m\bfd={\bf 0}$, where $\bfd=\bigl(v_{p_1}(d),\dots,v_{p_t}(d)\bigr)^\rmT$.
\end{itemize}
R\'edei showed that
\[\fct{\theta: V\cap \bfN F}{C(F)[2]\cap C(F)^2}{d}{\cl[\fd]}\]
is a two-to-one onto homomorphism, where $(d)=\fd^2$.
Thus
\begin{equation}\label{eq:Redei_rank}
h_4(m)=t-1-\rank(\bfR_m).
\end{equation}
See \cite{Redei1934} and \cite[Example~2.6]{LiYu2020}.

\begin{example}
Let $n=p_1\cdots p_k$ be an odd positive square-free integer.
When $n\equiv 1\bmod 4$, we have
\begin{align*}
\bfR_n&=\bfA+\bfD_{-1},&
\bfR_{-n}&=\begin{pmatrix}
\bfA&\bfb_2\\
\bfb_{-1}^\rmT&\aleg{2}{n}
\end{pmatrix},\\
\bfR_{2n}&=\begin{pmatrix}
\bfA+\bfD_{-2}&\bfb_2\\
\bfb_{2}^\rmT&\aleg{2}{n}
\end{pmatrix},&
\bfR_{-2n}&=\begin{pmatrix}
\bfA+\bfD_{2}&\bfb_2\\
\bfb_{-2}^\rmT&\aleg{2}{n}
\end{pmatrix}.
\end{align*}
When $n\equiv -1\bmod 4$, we have
\begin{align*}
\bfR_n&=\begin{pmatrix}
\bfA+\bfD_{-1}&\bfb_2\\
\bfb_{-1}^\rmT&\aleg{2}{n}
\end{pmatrix},&
\bfR_{-n}&=\bfA,\\
\bfR_{2n}&=\begin{pmatrix}
\bfA+\bfD_{-2}&\bfb_2\\
\bfb_{-2}^\rmT&\aleg{2}{n}
\end{pmatrix},&
\bfR_{-2n}&=\begin{pmatrix}
\bfA+\bfD_{2}&\bfb_2\\
\bfb_{2}^\rmT&\aleg{2}{n}
\end{pmatrix}.
\end{align*}
\end{example}

\subsection{Tame kernel}
Denote by $K_2\CO_F$ is the tame kernel of $F$.
We list the results of $r_4(K_2\CO_F)$ that we will use.
Assume that $|m|>2$.

\begin{theorem}[{\cite{BrowkinSchinzel1982}}]
The subgroup $K_2\CO_F[2]$ is generated by the Steinberg symbols
\begin{itemize}
\item $\set{-1,d}, d\mid m$;
\item $\{-1,u+\sqrt{m}\}$, where $m=u^2-cw^2$ for some $c=-1,\pm2$ and $u,w\in\BN$.
\end{itemize}
Denote by $k$ the number of odd prime factors of $m$.
If $m>2$, then 
\[r_2(K_2\CO_F)=k+\log_2\#\bigl(\set{\pm1,\pm2}\cap\bfN F\bigr).\]
If $m<-2$, then 
\[r_2(K_2\CO_F)=k-1+\log_2\#\bigl(\set{1,2}\cap\bfN F\bigr).\]
\end{theorem}

Denote by $m'=|m|/(2,m)$ the positive odd part of $m$.
If $2\notin\bfN F$, set $V_2=\emptyset$. 
\begin{theorem}[{\cite[Theorem~3.4]{Qin1995b}}]\label{thm:K_2_odd}
Suppose that $m>2$. 
Denote by $V_1$ the set of positive $d\mid m'$ satisfying: there exists $\varepsilon=\pm1$ or $\pm2$ such that $(d,-m)_p=\leg{\varepsilon}{p}, \forall p\mid m'$.
If $2\in\bfN F$, then write $m=2\mu^2-\lambda^2, \mu,\lambda\in\BN$ and denote by $V_2$ the set of positive $d\mid m'$ satisfying: there exists $\varepsilon=\pm1$ such that $(d,-m)_p=\leg{\varepsilon \mu}{p}, \forall p\mid m'$. We have
\[2^{r_4(K_2\CO_F)+1}=\#V_1+\#V_2.\]
\end{theorem}

\begin{theorem}[{\cite[Theorem~4.1]{Qin1995a}}]\label{thm:K_2_even}
Suppose that $m<-2$. 
Denote by $V_1$ the set of $d\mid m'$ satisfying: there exists $\varepsilon=1$ or $2$ such that $(d,-m)_p=\leg{\varepsilon}{p}, \forall p\mid m'$.
If $2\in\bfN F$, then write $m=2\mu^2-\lambda^2, \mu,\lambda\in\BN$ and denote by $V_2$ the set of $d\mid m'$ satisfying: $(d,-m)_p=\leg{\mu}{p}, \forall p\mid m'$. We have
\[2^{r_4(K_2\CO_F)+2}=\#V_1+\#V_2.\]
\end{theorem}

Denote by $\bfB=\bfB_m=\bfA_{m'}+\bfD_{m/m'}$, where $\bfA_{m'}$ is defined as \eqref{eq:monsky}.
If $m>0$, then
\begin{equation}\label{eq:v1v2_positive}
V_1=\set{\psi(\bfd): \bfB\bfd=\bfb_{\pm1},\bfb_{\pm2}},\qquad
V_2=\set{\psi(\bfd): \bfB\bfd=\bfb_{\pm\mu}}.
\end{equation}
If $m<0$, then
\begin{equation}\label{eq:v1v2_negative}
\begin{split}
V_1&=\set{\psi(\bfd): \bfB\bfd={\bf0},\bfb_2}\cup\set{-\psi(\bfd): \bfB\bfd=\bfb_{-1},\bfb_{-2}},\\
V_2&=\set{\psi(\bfd): \bfB\bfd=\bfb_{\mu}}\cup\set{-\psi(\bfd): \bfB\bfd=\bfb_{-\mu}}.
\end{split}
\end{equation}

\section{The odd case}

\begin{lemma}\label{lem:non-deg}
If $s_2(n)=2$, then $\rank_\BZ E_n(\BQ)=0,\Sha(E_n/\BQ)[2^\infty]\cong(\BZ/2\BZ)^2$ holds if and only if the Cassels pairing on the pure $2$-Selmer group $\Sel_2(E_n)/E_n(\BQ)[2]$ is non-degenerate.
\end{lemma}
\begin{proof}
See \cite[pp~2146, 2157]{Wang2016}.
\end{proof}

\begin{theorem}\label{thm:main1}
Let $n=p_1\cdots p_k\equiv 1\bmod 8$ be a square-free positive integer with $p_i\equiv \pm1\bmod8$.
The following are equivalent:
\begin{enumerate}[(i)]
\item $h_4(-n)=1, h_4(n)=0$ and $\leg{-\mu}{d}=-1$;
\item $\rank_\BZ E_n(\BQ)=0$ and $\Sha(E_n/\BQ)[2^\infty]\cong(\BZ/2\BZ)^2$.
\end{enumerate}
Here, $d\neq 1$ is the unique positive divisor of $n$ such that $(d,n)_v=1,\forall v$ and $n=2\mu^2-\tau^2$ where $\mu\equiv d\bmod4$.
\end{theorem}

\begin{proof}
In this case, $\bfD_2=\bfO$ and the Monsky matrix \eqref{eq:monsky_matrix1} is
\[\bfM_n=\begin{pmatrix}
\bfA&\\&\bfA+\bfD_{-1}
\end{pmatrix}.\]
Since $\bfA{\bf1}=(\bfA^\rmT+\bfD_{-1}){\bf1}={\bf0}$, we have $\rank(\bfA)\le k-1$ and $\rank(\bfA+\bfD_{-1})\le k-1$.
By \eqref{eq:monsky_rank}, we have
\[s_2(n)=2\iff \rank(\bfA)=\rank(\bfA+\bfD_{-1})=k-1.	\]
Since the R\'edei matrices
\[\bfR_{-n}=\left(\begin{smallmatrix}
\bfA&{\bf0}\\\bfb_{-1}^\rmT&0
\end{smallmatrix}\right),\quad
\bfR_n=\bfA+\bfD_{-1}\]
and note that ${\bf1}^\rmT\bfA=\bfb_{-1}^\rmT$, we have $\rank(\bfR_{-n})=\rank(\bfA)$ and
\[s_2(n)=2\iff h_4(-n)=1\text{ and }h_4(n)=0\]
by \eqref{eq:Redei_matrix} and \eqref{eq:Redei_rank}.

From the definition of $d$, we know that $\Ker(\bfA+\bfD_{-1})=\set{{\bf0},\psi^{-1}(d)}$.
The pure $2$-Selmer group of $E_n$ is
\[\set{(1,1,1),(1,n,n),(d,1,d),(d,n,nd)}\]
by \eqref{eq:monsky_isom1}.
Denote by $\Lambda=(1,n,n)$ and $\Lambda'=(d,1,d)$.
Then $D_\Lambda$ is defined by
\begin{equation}
\begin{cases}
H_1:& -t^2+u_2^2-u_3^2=0,\\
H_2:& -nt^2+nu_3^2-u_1^2=0,\\
H_3:& 2nt^2+u_1^2-nu_2^2=0.
\end{cases}
\end{equation}
Recall that $n=2\mu^2-\tau^2$.
Then $\mu$ is odd and $n=u^2-2w^2$, where $u=2\mu-\tau$ and $w=-\mu+\tau$.
Choose 
\begin{align*}
Q_1&=(0,1,1)\in H_1(\BQ),& L_1&=u_2-u_3,\\
Q_3&=(w,n,u)\in H_3(\BQ),& L_3&=2wt+u_1-uu_2.
\end{align*}
By Lemma~\ref{lem:cassels},
\[\pair{\Lambda,\Lambda'}=\prod_{p\mid 2n} \bigl(L_1L_3(P_p),d\bigr)_p\]
for any $P_p\in D_\Lambda(\BQ_p)$.

When $p\mid n$, take $(t,u_1,u_2,u_3)=(1,0,\sqrt{2},1)$ where $\sqrt{2}\equiv -u/w\bmod p$.
Then
\[L_1L_3(P_p)=(\sqrt{2}-1)(2w-\sqrt{2}u)\equiv 4(\sqrt{2}-1)w\equiv -4\mu\mod p\]
and $\bigl(L_1L_3(P_p),d\bigr)_p=(-\mu,d)_p$.
When $p=2$, take $(t,u_1,u_2,u_3)=(0,\sqrt{n},1,-1)$.
Then
\[\bigl(L_1L_3(P_2),d\bigr)_2=\bigl(2(\sqrt{n}-u),d\bigr)_2=(-\sqrt{n}-u,d)_2=(-\mu,d)_2\]
by Lemma~\ref{lem:2-symbol}.
Since $\mu\equiv d\bmod 4$, we have $(-\mu,d)_2=1$ and 
\[\pair{\Lambda,\Lambda'}=\prod_{p\mid d}(-\mu,d)_p=\leg{-\mu}{d}.\]
Conclude the result by Lemma~\ref{lem:non-deg}.
\end{proof}

\begin{lemma}\label{lem:2-symbol}
We have $(-u\pm \sqrt{n},-1)_2=(-\mu,-1)_2$.
\end{lemma}
\begin{proof}
See \cite[Lemma~3.1]{Qin1995b}.
Clearly, $(\mu,-1)_2=\pm 1$ if and only if $(-1,\pm \mu)_2=1$.
Thus the equation $X^2+Y^2=\pm \mu$ is solvable in $\BQ_2$ and so is $X^2+nY^2=\pm \mu$.
Let $x,y\in\BQ_2$ such that $x^2+ny^2=\pm \mu$.
Denote by
\[h=y,\quad g=\frac{x-wy}{\mu},\quad \alpha=g^2+h^2,\]
\[\theta=g^2-h^2+2gh,\quad \lambda=g^2-h^2-2gh,\]
\[\xi=\frac{g+h}{\alpha},\quad \eta=\frac{g-h}{\alpha},\]
\[x=-\xi+\lambda\eta,\ y=\alpha \xi,\ a=-\eta-\lambda\xi,\ b=\alpha \eta.\]
Then 
\[u \alpha+w \theta=\pm 1,\quad \theta^2+\lambda^2=2\alpha^2,
\quad \xi^2+\eta^2=\frac{2}{\alpha}\]
and
\[\begin{split}
&\left(\frac{x+y\sqrt{n}}{2}\right)^2+\left(\frac{a+b\sqrt{n}}{2}\right)^2\\
=&\frac14(\xi^2+\eta^2)(1+\lambda^2w^2+n\alpha^2-2\alpha\sqrt{n})\\
=&\frac{1}{2\alpha}\bigl((u\alpha+\theta w)^2+\lambda^2w^2+(u^2-2w^2)\alpha^2\bigr)-\sqrt{n}\\
=&u(u\alpha+\theta w)-\sqrt{n}=\pm u-\sqrt{n}.
\end{split}\]
Hence we have $(\mu,-1)_2=(u-\sqrt{n},-1)_2=(u+\sqrt{n},-1)_2$.
\end{proof}

%\begin{example}
%If $n=7\times 31=2\times 12^2-11^2$, then $d=7$, $\prod_{p\mid 2d}(\mu,d)_p=(13,7)_7(13,7)_2=-1$ and $n$ is non-congruent.
%\end{example}

When all prime factors of $n$ are $\equiv 1\bmod8$, we have the following corollary generalizing \cite[Lemma~4.2]{Qin2021}.
\begin{corollary}\label{cor:main1}
Let $n=p_1\cdots p_k\equiv 1\bmod 8$ be a square-free positive integer with $p_i\equiv 1\bmod8$.
The following are equivalent:
\begin{enumerate}[(i)]
\item $r_4(K_2\CO_{\BQ(\sqrt{n})})=0$;
\item $\rank_\BZ E_n(\BQ)=0$ and $\Sha(E_n/\BQ)[2^\infty]\cong(\BZ/2\BZ)^2$.
\end{enumerate}
\end{corollary}

\begin{proof}
Since all $p_i\equiv1\bmod8$, we have $\bfD_\varepsilon=\bfO$ and $\bfb_\varepsilon={\bf0}$ for any $\varepsilon=\pm1,\pm2$.
Note that $d=n$ in Theorem~\ref{thm:main1}, the condition (ii) is equivalent to $\rank(\bfA)=k-1$ and $\leg{\mu}{n}=\leg{-\mu}{n}=-1$.

Denote by $F=\BQ(\sqrt{n})$.
Combining with Theorem~\ref{thm:K_2_odd} and \eqref{eq:v1v2_positive}, we have
\[2^{r_4(K_2\CO_F)+1}=\#\set{\bfd: \bfA\bfd={\bf0}}+\#\set{\bfd:\bfA\bfd=\bfb_{|\mu|}=\bfb_\mu}\]
and $\Ker\bfA\supseteq\set{{\bf0},{\bf1}}$.
This implies that $r_4(K_2\CO_F)=0$ if and only if $\Ker\bfA=\set{{\bf0},{\bf1}}$ and $\bfA\bfd=\bfb_{\mu}$ is not solvable. Once we have $\Ker\bfA=\set{{\bf0},{\bf1}}$, we will get $\rank(\bfA)=k-1$. Since ${\bf1}^\rmT\bfA={\bf0}^\rmT$, we will have $\Im\bfA=\set{\bfd:{\bf1}^\rmT\bfd=0}$.
Hence $\bfA\bfd=\bfb_{\mu}$ is solvable if and only if $\leg{\mu}{n}=1$.
\end{proof}

One can propose many equivalent conditions for (i), which generalizes \cite{BarrucandCohn1969} and \cite[Lemma-Definition~1]{LiTian2000}. See also \cite[Theorem~4.2]{Wang2016}.

\begin{proposition}
Let $n$ be a square-free positive integer with prime factors congruent to $1$ modulo $8$. If $h_4(-n)=0$, then the following are equivalent:
\begin{enumerate}[(i)]
\item $2\mid b$, where $n=a^2+8b^2$;
\item $\leg{1+\sqrt{2}}{n}=1$;
\item $\leg{1+\sqrt{-1}}{n}=1$;
\item $\leg{\sqrt{2}}{n}=(-1)^{\frac{n-1}{8}}$;
\item $|u|\equiv 1\bmod 4$;
\item $|\mu|\equiv 1\bmod 4$;
\item $h_8(-n)=0$;
\item $\leg{\mu}{n}=1$.
\end{enumerate}
\end{proposition}
\begin{proof}
For any $p\mid n$, we have $\leg{a}{p}=\leg{\sqrt{2}(1+\sqrt{-1})^2b}{p}=\leg{\sqrt{2}b}{p}$. Thus $\leg{a}{n}=\leg{\sqrt{2}b}{n}$.
Since
\[\leg{a}{n}=\leg{n}{|a|}=\leg{a^2+8b^2}{|a|}=\leg{2}{|a|}=(-1)^{\frac{a^2-1}8}\]
and
\[\leg{b}{n}=\leg{n}{b'}=\leg{a^2+8b^2}{b'}=1,\quad b'=\text{odd part of }|b|,\]
we have $(-1)^{\frac{a^2-1}8}=\leg{\sqrt{2}}{n}$.
By \cite[Lemma-Definition~1]{LiTian2000}, $2\mid b$ if and only if
\[\begin{split}
1&=(-1)^{\frac{n-a^2}{8}}
=(-1)^{\frac{n-1}{8}}\leg{\sqrt{2}}{n}\\
&=\prod_{p\mid n} (-1)^{\frac{p-1}{8}}\leg{\sqrt{2}}{p}\\
&=\prod_{p\mid n}\leg{1+\sqrt{2}}{p}=\leg{1+\sqrt{2}}{n}.
\end{split}\]
This is equivalent to (iii)-(vii) by \cite[Theorem~4.2]{Wang2016}.
Note that
	\[\leg{\mu}{n}=\leg{n}{|\mu|}=\leg{2\mu^2-\tau^2}{|\mu|}=\leg{-1}{|\mu|},\]
we have (vi)$\iff$(viii).
\end{proof}

\section{The even case}

\begin{theorem}\label{thm:main2}
Let $n=2p_1\cdots p_k\equiv 2\bmod 8$ be a square-free positive integer with $p_i\equiv \pm1\bmod8$.
The following are equivalent:
\begin{enumerate}[(i)]
\item $h_4(-n/2)=1$ and $\leg{2-\sqrt{2}}{|d|}=-1$;
\item $\rank_\BZ E_n(\BQ)=0$ and $\Sha(E_n/\BQ)[2^\infty]\cong(\BZ/2\BZ)^2$.
\end{enumerate}
Here, $d\neq 1$ is the unique divisor of $n$ such that $(d,n)_v=1,\forall v$ and $d\equiv1\bmod8$.
\end{theorem}

\begin{proof}
In this case, $\bfD_2={\bfO}$ and the Monsky matrix \eqref{eq:monsky_matrix2} is
\[\bfM_n=\begin{pmatrix}
	\bfA^\rmT&\bfD_{-1}\\
	&\bfA
\end{pmatrix}.\]
Since $\bfA{\bf1}={\bf0}$, we have $\rank(\bfA)\le k-1$.
The equation $\bfM_n\left(\begin{smallmatrix}
\bfx\\\bfy
\end{smallmatrix}\right)={\bf0}$ can be written as 
\[\bfA^\rmT \bfx=\bfD_{-1}\bfy,\quad \bfA\bfy={\bf0}.\]
If $\bfy={\bf0}$, then $\bfA^\rmT\bfx={\bf0}$, which has at least two solutions.
If $\bfy={\bf1}$, then $\bfA^\rmT(\bfx+{\bf1})={\bf0}$, which has at least two solutions. Hence $s_2(n)=\dim_{\BF_2}\Ker\bfM_n\ge 2$ and we have
\[s_2(n)=2\iff \rank(\bfA)=k-1.\]
Denote by $n'=n/2$.
Since the R\'edei matrix
\[\bfR_{-n'}=\left(\begin{smallmatrix}
\bfA&{\bf0}\\\bfb_{-1}^\rmT&0
\end{smallmatrix}\right)\]
and note that ${\bf1}^\rmT\bfA=\bfb_{-1}^\rmT$, we have
\[s_2(n)=2\iff h_4(-n/2)=1	\]
by \eqref{eq:Redei_matrix} and \eqref{eq:Redei_rank}.

From the definition of $d$, we know that $\Ker\bfA^\rmT=\set{{\bf0},\psi^{-1}(|d|)}$.
The pure $2$-Selmer group of $E$ is
\[\set{(1,1,1),(d,1,d),(d,n',dn'),(1,n',n')}\]
by \eqref{eq:monsky_isom2}.
Denote by $\Lambda=(1,n',n')$ and $\Lambda'=(d,1,d)$.
Then $D_\Lambda$ is defined by
\begin{equation}
\begin{cases}
H_1:& -2t^2+u_2^2-u_3^2=0,\\
H_2:& -2n' t^2+n' u_3^2-u_1^2=0,\\
H_3:& 4n' t^2+u_1^2-n' u_2^2=0.
\end{cases}
\end{equation}
Choose 
\begin{align*}
Q_1&=(0,1,1)\in H_1(\BQ),& L_1&=u_2-u_3,\\
Q_3&=(n'-1,4n',2n'+2)\in H_3(\BQ),& L_3&=2(n'-1)t+2u_1-(n'+1)u_2.
\end{align*}
By Lemma~\ref{lem:cassels},
\[\pair{\Lambda,\Lambda'}=\prod_{p\mid 2n} \bigl(L_1L_3(P_p),d\bigr)_p\]
for any $P_p\in D_\Lambda(\BQ_p)$.

When $p\mid n$, take $(t,u_1,u_2,u_3)=(1,0,2,\sqrt{2})$.
Then
\[\bigl(L_1L_3(P_p),d\bigr)_p=\bigl(-2(2-\sqrt{2}),d\bigr)_p=(\sqrt{2}-2,d)_p.\]
When $p=2$, take $(t,u_1,u_2,u_3)=(0,\sqrt{n'},1,-1)$.
Then
\[L_1L_3(P_2)=2(2\sqrt{n'}-n'-1)=-2(\sqrt{n'}-1)^2\]
and $\bigl(L_1L_3(P_2),d\bigr)_2=(-2,d)_2=(-1,d)_2$.
Hence
\[\pair{\Lambda,\Lambda'}=(-1,d)_2\prod_{p\mid d}(\sqrt{2}-2,d)_p=\leg{2-\sqrt{2}}{|d|}.\]
Conclude the result by Lemma~\ref{lem:non-deg}.
\end{proof}

\begin{corollary}\label{cor:main2}
Let $n=2p_1\cdots p_k\equiv 2\bmod 8$ be a square-free positive integer with $p_i\equiv 1\bmod8$.
The following are equivalent:
\begin{enumerate}[(i)]
\item $r_4(K_2\CO_{\BQ(\sqrt{-n})})=0$;
\item $\rank_\BZ E_n(\BQ)=0$ and $\Sha(E_n/\BQ)[2^\infty]\cong(\BZ/2\BZ)^2$.
\end{enumerate}
\end{corollary}

\begin{proof}
Since all $p_i\equiv1\bmod8$, we have $\bfD_\varepsilon=\bfO$ and $\bfb_\varepsilon={\bf0}$ for any $\varepsilon=\pm1,\pm2$.
Note that $d=n'$ in Theorem~\ref{thm:main2},
the condition (ii) is equivalent to $\rank(\bfA)=k-1$ and $\leg{2-\sqrt{2}}{n'}=-1$.

Denote by $F=\BQ(\sqrt{-n})$ and write $-n=2\mu^2-\lambda^2,\mu,\lambda\in\BN$.
Combining with Theorem~\ref{thm:K_2_even} and \eqref{eq:v1v2_negative}, we have
\[2^{r_4(K_2\CO_F)+2}=2\#\set{\bfd: \bfA\bfd={\bf0}}+2\#\set{\bfA\bfd=\bfb_{|\mu|}=\bfb_\mu}\]
and $\Ker\bfA\supseteq\set{{\bf0},{\bf1}}$.
Therefore, $r_4(K_2\CO_F)=0$ if and only if $\Ker\bfA=\set{{\bf0},{\bf1}}$ and $\bfA\bfd=\bfb_{\mu}$ is not solvable. Once we have $\Ker\bfA=\set{{\bf0},{\bf1}}$, we will get $\rank(\bfA)=k-1$. Since ${\bf1}^\rmT\bfA={\bf0}^\rmT$, we will have $\Im\bfA=\set{\bfd:{\bf1}^\rmT\bfd=0}$.
Hence $\bfA\bfd=\bfb_{\mu}$ is solvable if and only if $\leg{\mu}{n'}=1$.

Denote by $u=\lambda-\mu$ and $w=-\lambda/2+\mu$.
Then $n'=u^2-2w^2$ and
\[\leg{\mu}{n'}=\leg{u+2w}{n'}=\leg{(2\pm\sqrt{2})w}{n'}.\]
The result then follows from
\[\leg{w}{n'}=\leg{n'}{w'}=\leg{u^2-2w^2}{w'}=\leg{u^2}{w'}=1,\]
where $w'$ is the positive odd part of $w$.
\end{proof}

\end{document}